\documentclass[12pt]{amsart}

\usepackage{amscd,latexsym,amsthm,amsfonts,amssymb,amsmath,amsxtra}
\usepackage[mathscr]{eucal}
\usepackage{hyperref}
\renewcommand\theequation{\thesection.\arabic{equation}}

\newcommand{\CH}{{\mathcal {H}}}

\newcommand{\RG}{{\mathrm {G}}}
\newcommand{\RH}{{\mathrm {H}}}

\newcommand{\GL}{{\mathrm{GL}}}

\newcommand{\SO}{{\mathrm{SO}}}

\newcommand{\Sp}{{\mathrm{Sp}}}

\newcommand{\wt}{\widetilde}

\newcommand{\bs}{\backslash}

\newtheorem{thm}{Theorem}[section]
\newtheorem{cor}[thm]{Corollary}

\newtheorem {ques/conj}[thm]{Question/Conjecture}

\makeatletter

\newcommand{\Rmnum}[1]{\expandafter\@slowromancap\romannumeral #1@}
\makeatother

\begin{document}
\renewcommand{\theequation}{\arabic{equation}}
\numberwithin{equation}{section}

\title{Model transition under local theta correspondence}

\author{Baiying Liu}
\address{Department of Mathematics\\
University of Utah\\
155 S 1400 E Room 233, Salt Lake City, UT 84112-0090, USA.}
\email{liu@math.utah.edu}


\date{\today}


\thanks{The author is supported in part by NSF Grants DMS-1302122, and in part by a postdoc research fund from Department of Mathematics, University of Utah}

\maketitle

\begin{abstract}
We study model transition for representations occurring in the local theta correspondence between split even special orthogonal groups and symplectic groups, over a non-archimedean local field of characteristic zero. 
\end{abstract}


\section{Introduction}

Let $F$ be a non-archimedean local field of characteristic zero.
Let $\RG$ be a reductive group defined over $F$. An important way to characterize irreducible admissible representations of $\RG(F)$ is to consider various kinds of models. For example, non-degenerate Whittaker models of generic representations play significant roles in the theory of local factors of representations of $\RG(F)$ and in the theory of automorphic forms. Given a reductive dual pair $(\RG, \RH)$, and consider the local theta correspondence between representations of $\RG(F)$ and of $\RH(F)$, an interesting question is that how models of representations of $\RG(F)$ and of $\RH(F)$ are related under local theta correspondence. 
 
In \cite{JNQ10a}, Jiang, Nien and Qin proved that under the local theta correspondence, for certain representations, the generalized Shalika
model on $\SO_{4n}(F)$ is corresponding to the symplectic linear model on $\Sp_{4n}(F)$, and conjectured that it is true for general representations (see \cite[p. 542]{JNQ10b}). In this paper, 
we prove some results related to this conjecture (see Theorem \ref{main}).
More precisely, we introduce generalized symplectic linear models on $\Sp_{4m}(F)$, which generalizes the symplectic linear models, and study the relations between the generalized Shalika
models on $\SO_{4n}(F)$ and the generalized symplectic linear models on $\Sp_{4m}(F)$ under the local theta correspondence. 
A special case ($m=n$) of this result is proved by Hanzer (\cite{H15}) independently. 
We also introduce generalized Shalika models on $\Sp_{4n}(F)$ and generalized orthogonal linear models on $\SO_{4m}(F)$, and study the relations between them under the local theta correspondence. 

M{\oe}glin (\cite{M98}), Gomez and Zhu (\cite{GZ14}) studied the local theta lifting of generalized Whittaker models associated to nilpotent orbits.
Note that the generalized Shalika models are indeed generalized Whittaker models associated to certain nilpotent orbits, but the generalized symplectic/orthogonal linear models are not. For example, the generalized Shalika models on $\SO_{4n}(F)$ are  generalized Whittaker models associated to the nilpotent orbits parametrized by the partition $[2^{2n}]$. By \cite{M98} and \cite{GZ14}, the full local theta lift on $\Sp_{2k}(F)$ has a nonzero generalized Whittaker model associated to a nilpotent orbit parametrized by the partition $[3^{2n}1^{\ell}]$, if 
$6n+\ell=2k$.
In general it is not known whether the small theta lift on $\Sp_{2k}(F)$ (if nonzero) would also carry this model. 

This paper is organized as follows. 
In Section 2, we give the definitions for various models for representations of split even special orthogonal groups and symplectic groups, and introduce the main result Theorem \ref{main}. In Section 3 and 4, we prove Part (1) and Part (2) of Theorem \ref{main} respectively. In Section 5, we consider the converse of Theorem \ref{main}, and discuss some related results.

\textbf{Acknowledgments.} The main results of this paper were worked out when the author was a graduate student at University of Minnesota. The author would like to thank his advisor Prof. Dihua Jiang for constant support, encouragement and helpful discussions. The author also would like to thank Prof. Gordan Savin, Prof. David Soudry and Prof. Chen-bo Zhu for helpful conversations, and Raul Gomez for helpful communications. The author also would like to thank the referee for helpful comments and suggestions on this paper.

\section{Models of representations}

In this section, we define various models for representations of split even special orthogonal groups and symplectic groups. 
 
For any positive integer $k$, let $v_k$ be the $k\times k$ matrix with $1$'s in the second diagonal and zero's elsewhere. 
Let 
$$\SO_{2\ell} = \{g \in \GL_{2\ell} \ | \ {}^tg v_{2\ell} g =v_{2\ell}\}$$
be the split even special orthogonal group. And let $SO_{2\ell} = \SO_{2\ell}(F)$. Let 
$$\Sp_{2\ell} = \{g \in \GL_{2\ell} | {}^tg J_{2\ell} g =J_{2\ell}\}$$
be the symplectic group, where $J_{2\ell} =\begin{pmatrix}
        0 & v_{\ell} \\
       -v_{\ell} & 0 
\end{pmatrix}$. And let $Sp_{2\ell}=\Sp_{2\ell}(F)$. 

Let $P_{2n} = M_{2n} N_{2n}$ be the Siegel parabolic subgroup of $SO_{4n}$. Then elements in $P_{2n}$ have the following form
$$
(g,X) = m(g) n(X) = \left(\begin{array}{cc}
        g & 0 \\
        0 & g^* 
      \end{array}
\right) \left(\begin{array}{cc}
        I_n & X \\
        0 & I_n 
      \end{array}
\right)
$$
where $g \in \GL_{2n}(F)$, $g^*= v_{2n} {}^tg^{-1} v_{2n}$, and $X \in M_{(2n) \times (2n)}(F)$ satisfying ${}^tX = - v_{2n} X v_{2n}$.
The generalized Shalika group of $SO_{4n}$ is defined to be
$$\mathcal{H} = \mathcal{H}_{2n} = M'_{2n} N_{2n}= \{ (g,X) \in P_{2n} \ | \ g \in Sp_{2n} \}.$$ 
Define a character of $\CH$ as follows
\begin{align*}
\psi_{\mathcal{H}} ((g,X)) & = \psi(\frac{1}{2} tr(J_{2n} X v_{2n}))\\
& = \psi(\frac{1}{2} tr(\begin{pmatrix}
-I_n & 0 \\
0 & I_n
\end{pmatrix} X)).
\end{align*}
Then an irreducible admissible representation $(\sigma, V_{\sigma})$ of $SO_{4n}$ has a generalized Shalika model if 
$$\mathrm{Hom}_{SO_{4n}}(V_{\sigma}, \mathrm{Ind}^{SO_{4n}}_{\mathcal{H}} (\psi_{\mathcal{H}}))
 = \mathrm{Hom}_{\mathcal{H}}(V_{\sigma}, \psi_{\mathcal{H}}) \neq 0.$$ The nonzero elements in the 
$\mathrm{Hom}$ space are called generalized Shalika functionals or $\psi_{\mathcal{H}}$-functionals of $\sigma$. In 
\cite{N10}, Nien proved the uniqueness of the generalized Shalika models.

Let $Q_{2n} = L_{2n} V_{2n}$ be the Siegel parabolic subgroup of $Sp_{4n}$. Then elements in $Q_{2n}$ have the following form
$$
(g,X) = m(g) n(X) = \left(\begin{array}{cc}
        g & 0 \\
        0 & g^* 
      \end{array}
\right) \left(\begin{array}{cc}
        I_n & X \\
        0 & I_n 
      \end{array}
\right)
$$
where $g \in \GL_{2n}(F)$, $g^*= v_{2n} {}^tg^{-1} v_{2n}$, and $X \in M_{(2n) \times (2n)}(F)$ satisfying ${}^tX = v_{2n} X v_{2n}$.
The generalized Shalika group of $Sp_{4n}$ is defined to be
$$\wt{\mathcal{H}} = \wt{\mathcal{H}}_{2n} = L'_{2n} V_{2n}= \{ (g,X) \in P_{2n} \ | \ g \in SO_{2n} \}.$$ 
Let 
$$\psi_{\wt{\mathcal{H}}} ((g,X)) = \psi(\frac{1}{2} tr(X))$$
be a character of $\wt{\mathcal{H}}$. An irreducible admissible representation $(\pi, V_{\pi})$ of 
$Sp_{4n}$ has a generalized Shalika model if 
$$\mathrm{Hom}_{Sp_{4n}}(V_{\pi}, \mathrm{Ind}^{Sp_{4n}}_{\wt{\mathcal{H}}} (\psi_{\wt{\mathcal{H}}}))
 = \mathrm{Hom}_{\wt{\mathcal{H}}}(V_{\pi}, \psi_{\wt{\mathcal{H}}}) \neq 0.$$ The nonzero elements in the 
$\mathrm{Hom}$ space are called generalized Shalika functionals or $\psi_{\wt{\mathcal{H}}}$-functionals of $\pi$. 

Consider $SO_{2n}\times SO_{4m-2n}$ as a subgroup of $SO_{4m}$ ($2m > n$) via the embedding 
$$(\begin{pmatrix}
A_1 & B_1 \\
C_1 & D_1
\end{pmatrix}, \begin{pmatrix}
A_2 & B_2 \\
C_2 & D_2
\end{pmatrix}) \mapsto \begin{pmatrix}
A_1 & 0 & 0 & B_1\\
0 & A_2 & B_2 & 0\\
0 & C_2 & D_2 & 0\\
C_1 & 0 & 0 & D_2
\end{pmatrix}.$$
An irreducible admissible representation $(\sigma, V_{\sigma})$ of $SO_{4m}$ ($2m \geq n$) has a generalized 
orthogonal linear model if
\begin{align*}
& \ \mathrm{Hom}_{\SO_{4m}(F)}(V_{\sigma}, \mathrm{Ind}^{\SO_{4m}(F)}_{\SO_{2n}(F) \times \SO_{4m-2n}(F)} (1)\\
 = & \ \mathrm{Hom}_{\SO_{2n}(F) \times \SO_{4m-2n}(F)}(V_{\sigma}, 1) \\
 \neq & \ 0.
\end{align*}
The nonzero elements in the 
$\mathrm{Hom}$ space are called generalized orthogonal linear functionals of $\pi$.

Similarly, consider $Sp_{2n}\times Sp_{4m-2n}$ as a subgroup of $Sp_{4m}$ ($2m > n$) via the embedding 
$$(\begin{pmatrix}
A_1 & B_1 \\
C_1 & D_1
\end{pmatrix}, \begin{pmatrix}
A_2 & B_2 \\
C_2 & D_2
\end{pmatrix}) \mapsto \begin{pmatrix}
A_1 & 0 & 0 & B_1\\
0 & A_2 & B_2 & 0\\
0 & C_2 & D_2 & 0\\
C_1 & 0 & 0 & D_2
\end{pmatrix}.$$
An irreducible admissible representation $(\pi, V_{\pi})$ of $Sp_{4m}$ ($2m \geq n$) has a generalized 
symplectic linear model if
\begin{align*}
& \mathrm{Hom}_{Sp_{4m}}(V_{\pi}, \mathrm{Ind}^{Sp_{4m}}_{Sp_{2n}\times Sp_{4m-2n}} (1))\\
 = \ & \mathrm{Hom}_{Sp_{2n} \times Sp_{4m-2n}}(V_{\pi}, 1)\\
  \neq \ & 0
\end{align*}
The nonzero elements in the
$\mathrm{Hom}$ space are called generalized symplectic linear functionals of $\pi$. 
When
$m = n$, the generalized symplectic linear model is usually called the symplectic linear model. 
In \cite{Zh10}, Zhang proves the uniqueness
of the symplectic linear models for the cases of $n=1,2$. The uniqueness of generalized symplectic linear models on $Sp_{4m}$ is unknown in general.

Let $\omega_{\psi}$ be the Weil representation corresponding to the reductive dual pair $(O_{2k}, Sp_{2\ell})$ with $O_{2k}$ being split. 
Following is the main result of this paper.

\begin{thm}\label{main}
\begin{enumerate}
\item Assume that $\sigma$ is an irreducible admissible representation of $SO_{4n}$ with a nonzero generalized Shalika model,
and $\pi$ is an irreducible admissible representation of $Sp_{4m}$ $(2m \geq n)$ which is corresponding to $\sigma$ under the 
local theta correspondence, that is,  
$$\mathrm{Hom}_{SO_{4n} \times Sp_{4m}} (\omega_{\psi}, \sigma \otimes \pi) \neq 0.$$ Then $\pi$ has a nonzero generalized 
symplectic linear model.
\item Assume that $\pi$ is an irreducible admissible representation of $Sp_{4n}$ with a nonzero generalized Shalika model,
and $\sigma$ is an irreducible admissible representation of $SO_{4m}$ $(2m \geq n)$ which is corresponding to $\pi$ under the local theta correspondence, that is,  
$$\mathrm{Hom}_{SO_{4m} \times Sp_{4n}} (\omega_{\psi}, \sigma \otimes \pi) \neq 0.$$
Then, $\sigma$ has a nonzero generalized orthogonal linear model.
\end{enumerate}
\end{thm}

In \cite{H15}, Hanzer proved the case of $m=n$ for Theorem \ref{main}, Part (1), independently. Some results related to the converse of Theorem \ref{main} will be discussed in Section 5. 

M{\oe}glin (\cite{M98}), Gomez and Zhu (\cite{GZ14}) considered the local theta lifting of generalized Whittaker models associated to nilpotent orbits. It turns out a nonzero generalized Shalika model for an irreducible admissible representation $\sigma$ of $SO_{4n}$ is indeed a generalized Whittaker model associated to a nilpotent orbit parametrized by the partition $[2^{2n}]$. Then they showed that the full local theta lift of $\sigma$ on $Sp_{2k}$ has a nonzero generalized Whittaker model associated to a nilpotent orbit parametrized by the partition $[3^{2n}1^{\ell}]$, if $6n+\ell=2k$.
Note that this later model is not a generalized symplectic linear model. And, in general it is not known that whether the small theta lift would also carry this model. 

\section{Proof of Theorem \ref{main}, Part (1)}

In this section, we prove Theorem \ref{main}, Part (1), using a method similar to that in \cite[Section 2.1]{JS03}.

Assume that $\sigma$ is an irreducible admissible representation of $SO_{4n}$ with a nonzero generalized Shalika model,
and $\pi$ is an irreducible admissible representation of $Sp_{4m}$ $(2m \geq n)$ which is corresponding to $\sigma$ under the 
local theta correspondence, i.e., 
$$\mathrm{Hom}_{SO_{4n} \times Sp_{4m}} (\omega_{\psi}, \sigma \otimes \pi) \neq 0.$$

Since $\sigma$ has a nonzero generalized Shalika model, i.e., 
$\mathrm{Hom}_{\mathcal{H}}(V_{\sigma}, \psi_{\mathcal{H}}) \neq 0,$ we have that
$$\mathrm{Hom}_{\mathcal{H} \times \Sp_{4m}(F)} (\omega_{\psi}, \psi_{\mathcal{H}} \otimes \pi) \neq 0.$$

Let $V$ be a $4n$-dimensional vector space over $F$, with the nondegenerate symmetric from $v_{4n}$. Fix 
a basis $\{e_1, \ldots, e_{2n}, e_{-2n}, \ldots, e_{-1}\}$ of $V$ over $F$, such that $(e_i, e_j) = (e_{-i}, e_{-j}) = 0,
(e_i, e_{-j}) = \delta_{ij}$, for $i, j = 1, \ldots, 2n$. Let $V^+ = \mathrm{Span}_F\{e_1, \ldots, e_{2n}\}$,
$V^- = \mathrm{Span}_F\{e_{-1}, \ldots, e_{-2n}\}$, then $V= V^+ + V^-$ is a complete polarization of $V$.
Let $W$ be a $4m$-dimensional symplectic vector space over $F$, with the symplectic from $J_{4m} = \left(\begin{array}{cc}
        0 & v_{2m} \\
       -v_{2m} & 0 
      \end{array}
\right)$. Fix 
a basis 
$$\{f_1, \ldots, f_{2m}, f_{-2m}, \ldots, f_{-1}\}$$
of $W$ over $F$, such that $(f_i, f_j) = (f_{-i}, f_{-j}) = 0,
(f_i, f_{-j}) = \delta_{ij}$, for $i, j = 1, \ldots, 2m$. Let $W^+ = \mathrm{Span}_F\{f_1, \ldots, f_{2m}\}$,
$W^- = \mathrm{Span}_F\{f_{-1}, \ldots, f_{-2m}\}$, then $W= W^+ + W^-$ is a complete polarization of $W$.
Then, we realize the Weil representation $\omega_{\psi}$ in the space $\mathcal{S}(V^- \otimes W) \cong \mathcal{S}(W^{2n})$.

First, we compute the Jacquet module $J_{N_{2n}, \psi_{\mathcal{H}}|_{N_{2n}}} ( \mathcal{S} (W^{2n}))$.
By \cite[Section 2.7]{MVW87},
for $\phi \in \mathcal{S}(W^{2n})$,
\begin{equation} \label{equ1}
\omega_{\psi} (n(X), 1) \phi (y_1, \ldots, y_{2n}) =  \psi(\frac{1}{2} tr(Gr(y_1, \ldots, y_{2n})v_{2n}X))  
\phi (y_1, \ldots, y_{2n}),
\end{equation}
where $Gr(y_1, \ldots, y_{2n})$ is the Gram matrix of vectors $y_1, \ldots, y_{2n}$. 
Since ${}^tX = - v_{2n} X v_{2n}$, 
\begin{align} \label{equ2}
\begin{split}
& \ \psi(\frac{1}{2} tr(Gr(y_1, \ldots, y_{2n})v_{2n}X)) \\
= & \ \psi(x_{11} (y_1, y_{2n}) + \cdots + x_{nn} (y_n, y_{n+1}) + x_{12} (y_2, y_{2n}) + \cdots),
\end{split}
\end{align}
where $\frac{1}{2} tr(Gr(y_1, \ldots, y_{2n})v_{2n}X)$ involves all $x_{ij}$ terms except $x_{1, 2n}$, $x_{2, 2n-1}$, $\ldots, x_{2n, 1}$, 
since they are zeros.
And, 
\begin{equation} \label{equ3}
\psi_{\mathcal{H}}(n(X)) = \psi(\frac{1}{2} tr(\left(\begin{array}{cc}
       -I_n & 0 \\
       0 & I_n 
      \end{array}
\right) X)) = \psi(-(x_{11} + \cdots + x_{nn})).
\end{equation}

By comparing \eqref{equ2}, \eqref{equ3}, let 
$$C_0 = \{ (y_1, \ldots, y_{2n}) \in W^{2n} | Gr(y_1, \ldots, y_{2n}) = -J_{2n}\}.$$
Then, by \eqref{equ1}, \eqref{equ2} and \eqref{equ3}, we claim that
\begin{equation} \label{claim1}
J_{N_{2n}, \psi_{\mathcal{H}}|_{N_{2n}}} ( \mathcal{S} (W^{2n})) \cong \mathcal{S}(C_0).
\end{equation}

Indeed, let $C = W^{2n} \bs C_0$, which is an open set, since $C_0$ is closed in $W^{2n}$. By \cite[Section 1.8]{BZ76}, there is an exact sequence
$$0 \rightarrow \mathcal{S}(C) \xrightarrow{i} \mathcal{S}(W^{2n}) \xrightarrow{r} \mathcal{S}(C_0) \rightarrow 0,$$
where $i$ is the canonical embedding that extends any function supported in $C$ by zero to the whole space $W^{2n}$,
and $r$ is the restriction map to $C_0$. By the exactness of Jacquet functors, there is also an exact sequence
$$0 \rightarrow J_{N_{2n}, \psi_{\mathcal{H}}|_{N_{2n}}} (\mathcal{S}(C)) \xrightarrow{i^*} 
J_{N_{2n}, \psi_{\mathcal{H}}|_{N_{2n}}}(\mathcal{S}(W^{2n})) \xrightarrow{r^*} 
J_{N_{2n}, \psi_{\mathcal{H}}|_{N_{2n}}}(\mathcal{S}(C_0)) \rightarrow 0.$$
By \eqref{equ1}, \eqref{equ2} and \eqref{equ3}, for any $(y_1, \ldots, y_{2n}) \in C$, the character
$$n(X) \mapsto \psi(\frac{1}{2} tr(Gr(y_1, \ldots, y_{2n})v_{2n}X)) - \psi_{\mathcal{H}}(n(X))$$
is nontrivial. Hence, for $\phi \in \mathcal{S}(C)$, there exists a large enough compact subgroup
$N_{\phi}$ of $N_{2n}$, such that 
$$\int_{N_{\phi}} (\omega_{\psi} (n(X), 1) - \psi_{\mathcal{H}}(n(X))) \phi  (y_1, \ldots, y_{2n}) dX$$
\begin{align*}
& =  \int_{N_{\phi}} (\psi(\frac{1}{2} tr(Gr(y_1, \ldots, y_{2n})v_{2n}X))  - \psi_{\mathcal{H}}(n(X))) \phi  (y_1, \ldots, y_{2n}) dX\\
& = 0.
\end{align*}
So, $J_{N_{2n}, \psi_{\mathcal{H}}|_{N_{2n}}} (\mathcal{S}(C)) = 0$, and 
$$J_{N_{2n}, \psi_{\mathcal{H}}|_{N_{2n}}}(\mathcal{S}(W^{2n})) \cong 
J_{N_{2n}, \psi_{\mathcal{H}}|_{N_{2n}}}(\mathcal{S}(C_0)) \cong \mathcal{S}(C_0).$$
This proves the claim in \eqref{claim1}.

Therefore,
\begin{equation} \label{equ5}
\mathrm{Hom}_{M'_{2n} \times Sp_{4m}} (\mathcal{S}(C_0), 1 \otimes \pi) \neq 0.
\end{equation}
Note that $M'_{2n} = \{ m(g) \in M_{2n}| g \in \Sp_{2n}\}.$

By \cite[Section 2.7]{MVW87}, for $\phi \in \mathcal{S}(C_0)$, $m(g) \in M'_{2n}$, $h \in Sp_{4m}$,
\begin{equation*}
\omega_{\psi} (m(g), h) \phi (y_1, \ldots, y_{2n}) = \phi ((y_1 h, \ldots, y_{2n} h) v_{2n} g v_{2n}).
\end{equation*}
Note that by Witt Theorem, $M'_{2n} \times Sp_{4m}$ acts transitively on $C_0$, and 
$$(f_{-2m + n -1}, \ldots, f_{-2m}, f_{2m}, \ldots, f_{2m-n+1})$$ 
is a representative. Let $R$ be the stabilizer of this 
representative. Then
\begin{equation*}
R = \{ (m(g), \left(\begin{array}{ccc}
       a & 0 & b\\
       0 & g^{-1} & 0\\
       c & 0 & d 
      \end{array}
\right))\ | \ g \in Sp_{2n}, \left(\begin{array}{cc}
       a & b \\
       c & d 
      \end{array}
\right) \in Sp_{4m-2n} \}.
\end{equation*}
Hence, $\mathcal{S}(C_0)$ is isomorphic to the compactly induced representation
$$\mathrm{c-Ind}^{M'_{2n} \times Sp_{4m}}_R (1).$$

Therefore, by \eqref{equ5}, 
\begin{equation*}
\mathrm{Hom}_{M'_{2n} \times Sp_{4m}}
 (\mathrm{c-Ind}^{M'_{2n} \times Sp_{4m}}_R (1), 1 \otimes \pi) \neq 0,
\end{equation*}
that is, 
$$
\mathrm{Hom}_{M'_{2n} \times Sp_{4m}}
 (1 \otimes \wt{\pi}, \mathrm{Ind}^{M'_{2n} \times Sp_{4m}}_R (1)) \neq 0,$$
where $\wt{\pi}$ is the contragredient of $\pi$. 
Note that in general, given admissible representations $\sigma$ and $\tau$ of a connected reductive group $G$, one has
$$\mathrm{Hom}_{G}(\sigma, \tau) \cong \mathrm{Hom}_{G}(\wt{\tau}, \wt{\sigma}).$$
And note that
$$\wt{\mathrm{c-Ind}^{M'_{2n} \times Sp_{4m}}_R (1)}
\cong \mathrm{Ind}^{M'_{2n} \times Sp_{4m}}_R (1),$$
since both $M'_{2n} \times Sp_{4m}$ and $R$ are unimodular.

Then, by Frobenius Reciprocity, we have
\begin{equation*} 
\mathrm{Hom}_R (1 \otimes \wt{\pi}, 1) \neq 0,
\end{equation*}
that is, 
\begin{equation*}
\mathrm{Hom}_{\Sp_{2n}(F) \times \Sp_{4m-2n}(F)} (\wt{\pi}, 1) \neq 0,
\end{equation*}
which means that $\wt{\pi}$ has a nonzero generalized symplectic linear model. It follows from an argument as in the proof of \cite[Theorem 17]{GRS99a} that $\pi$ also has a nonzero
generalized symplectic linear model. Indeed, by \cite[pages 91-92]{MVW87}, $\wt{\pi} \cong \pi^{\delta}$, 
where $\delta = v_{4m}$, and $\pi^{\delta}(g) = \pi(\delta g \delta^{-1}) =\pi(\delta g \delta)$.
Hence, $\pi^{\delta}$ has a nonzero generalized symplectic linear model. 
Assume that $l$ is a nonzero generalized symplectic linear functional of $\pi^{\delta}$.
Then, for any $(g_1, g_2) \in Sp_{2n} \times Sp_{4m-2n}$, 
$v \in V_{\pi}$, 
$$l(\pi(g_1, g_2) v) = l(\pi^{\delta}(g_1, g_2)^{\delta} v) = l(\pi^{\delta}(v_{2n} g_1 v_{2n},v_{4m-2n} g_2 v_{4m-2n}) v)
= l(v),$$
noting that $(v_{2n} g_1 v_{2n},v_{4m-2n} g_2 v_{4m-2n}) \in Sp_{2n} \times Sp_{4m-2n}$.
So, $l$ is also a nonzero generalized symplectic linear functional of $\pi$. Therefore, $\pi$ has a nonzero
generalized symplectic linear model. 

This completes the proof of Theorem \ref{main}, Part (1).

\section{Proof of Theorem \ref{main}, Part (2)}

In this section, we prove Theorem \ref{main}, Part (2), using a similar argument as in Section 3.

Assume that $\pi$ is an irreducible admissible representation of $Sp_{4n}$ with a nonzero generalized Shalika model,
and $\sigma$ is an irreducible admissible representation of $SO_{4m}$ $(2m \geq n)$ which is corresponding to $\pi$ under the local theta correspondence, i.e., 
$$\mathrm{Hom}_{SO_{4m} \times Sp_{4n}} (\omega_{\psi}, \sigma \otimes \pi) \neq 0.$$

Since $\pi$ has a nonzero generalized Shalika model, that is, 
$$\mathrm{Hom}_{\wt{\CH}}(V_{\pi}, \psi_{\wt{\CH}}) \neq 0,$$
we have that
$$\mathrm{Hom}_{SO_{4m} \times \wt{\CH}} (\omega_{\psi}, \sigma \otimes \psi_{\wt{\CH}}) \neq 0.$$

Let $V$ be a $4m$-dimensional vector space over $F$, with the nondegenerate symmetric from $v_{4m}$. Fix 
a basis $\{e_1, \ldots, e_{2m}, e_{-2m}, \ldots, e_{-1}\}$ of $V$ over $F$, such that $(e_i, e_j) = (e_{-i}, e_{-j}) = 0,
(e_i, e_{-j}) = \delta_{ij}$, for $i, j = 1, \ldots, 2m$. Let $V^+ = \mathrm{Span}_F\{e_1, \ldots, e_{2m}\}$,
$V^- = \mathrm{Span}_F\{e_{-1}, \ldots, e_{-2m}\}$, then $V= V^+ + V^-$ is a complete polarization of $V$.
Let $W$ be a $4n$-dimensional symplectic vector space over $F$, with the symplectic from $J_{4n} = \left(\begin{array}{cc}
        0 & v_{2n} \\
       -v_{2n} & 0 
      \end{array}
\right)$. Fix 
a basis 
$$\{f_1, \ldots, f_{2n}, f_{-2n}, \ldots, f_{-1}\}$$
of $W$ over $F$, such that $(f_i, f_j) = (f_{-i}, f_{-j}) = 0,
(f_i, f_{-j}) = \delta_{ij}$, for $i, j = 1, \ldots, 2n$. Let $W^+ = \mathrm{Span}_F\{f_1, \ldots, f_{2n}\}$,
$W^- = \mathrm{Span}_F\{f_{-1}, \ldots, f_{-2n}\}$, then $W= W^+ + W^-$ is a complete polarization of $W$.
Realize the Weil representation $\omega_{\psi}$ in the space $\mathcal{S}(V \otimes W^+) \cong \mathcal{S}(V^{2n})$. 

First, as in Section 3, we need to compute the Jacquet module of the Weil representation $J_{N_{2n}, \psi_{\mathcal{H}}|_{N_{2n}}} ( \mathcal{S} (V^{2n}))$.

By \cite[Section 2.6]{MVW87},
for $\phi \in \mathcal{S}(V^{2n})$,
\begin{equation} \label{equ14}
\omega_{\psi} (n(X), 1) \phi (x_1, \ldots, x_{2n}) =  \psi(\frac{1}{2} tr(Gr(x_1, \ldots, x_{2n})X v_{2n}))  
\phi (x_1, \ldots, x_{2n}),
\end{equation}
where $Gr(x_1, \ldots, x_{2n})$ is the Gram matrix of vectors $x_1, \ldots, x_{2n}$.
Since ${}^tX = v_{2n} X v_{2n}$, 
\begin{align} \label{equ15}
& \ \psi(\frac{1}{2} tr(Gr(x_1, \ldots, x_{2n})X v_{2n})) \\
=  & \ \psi(x_{11} (x_1, x_{2n}) + \cdots + x_{nn} (x_n, x_{n+1}) + \frac{1}{2} x_{1,2n} (x_1, x_1) + \cdots),
\end{align}
where $\frac{1}{2} tr(Gr(x_1, \ldots, x_{2n})X v_{2n})$ involves all $x_{ij}$ terms.
And, 
\begin{equation} \label{equ16}
\psi_{\mathcal{H}}(n(X)) = \psi(\frac{1}{2} tr(X)) = \psi(x_{11} + \cdots + x_{nn}).
\end{equation}

By comparing \eqref{equ14}, \eqref{equ15}, let 
$$C_0 = \{ (x_1, \ldots, x_{2n}) \in W^{2n} | Gr(x_1, \ldots, x_{2n}) = v_{2n}\}.$$
Then, by \eqref{equ14}, \eqref{equ15} and \eqref{equ16}, a similar argument (details omitted) as in Section 3 shows that
\begin{equation*}
J_{N_{2n}, \psi_{\mathcal{H}}|_{N_{2n}}} ( \mathcal{S} (V^{2n})) \cong \mathcal{S}(C_0).
\end{equation*}

Therefore,
\begin{equation} \label{equ17}
\mathrm{Hom}_{SO_{4m} \times M'_{2n}} (\mathcal{S}(C_0), \sigma \otimes 1) \neq 0.
\end{equation}
Note that $M'_{2n} = \{ m(g) \in M_{2n}| g \in \SO_{2n}\}.$

By \cite[Section 2.6]{MVW87}, for $\phi \in \mathcal{S}(C_0)$, $m(g) \in M'_{2n}$, $h \in SO_{4m}$,
\begin{equation*}
\omega_{\psi} (h, m(g)) \phi (x_1, \ldots, x_{2n}) = \phi ((h^{-1}x_1, \ldots, h^{-1}x_{2n})g).
\end{equation*}
Note that by Witt Theorem, $SO_{4m} \times M'_{2n}$ acts transitively on $C_0$, and 
$$(e_{2m-n+1}, \ldots, e_{2m}, e_{-2m}, \ldots, e_{-2m+n-1})$$ 
is a representative. Let $R$ be the stabilizer of this 
representative. Then
\begin{equation*}
R = \{ \left(\begin{array}{ccc}
       a & 0 & b\\
       0 & g^{-1} & 0\\
       c & 0 & d 
      \end{array}
\right)) | g \in \SO_{2n}, \left(\begin{array}{cc}
       a & b \\
       c & d 
      \end{array}
\right) \in SO_{4m-2n} \}.
\end{equation*}
Hence, $\mathcal{S}(C_0)$ is isomorphic to the compactly induced representation
$$\mathrm{c-Ind}^{\SO_{4m} \times M'_{2n}(F)}_R (1).$$

Hence, by \eqref{equ17}, 
\begin{equation*}
\mathrm{Hom}_{SO_{4m} \times M'_{2n}}
 (\mathrm{c-Ind}^{SO_{4m} \times M'_{2n}}_R (1), \sigma \otimes 1) \neq 0,
\end{equation*}
that is,  
$$\mathrm{Hom}_{SO_{4m} \times M'_{2n}}
 (\wt{\sigma} \otimes 1, \mathrm{Ind}^{SO_{4m} \times M'_{2n}}_R (1)) \neq 0,$$
as in Section 3. 
By Frobenius Reciprocity, we have
\begin{equation*}
\mathrm{Hom}_R (\wt{\sigma} \otimes 1, 1) \neq 0,
\end{equation*}
that is, 
\begin{equation*}
\mathrm{Hom}_{SO_{2n} \times SO_{4m-2n}} (\wt{\sigma}, 1) \neq 0.
\end{equation*}
So, $\wt{\sigma}$ has a nonzero generalized orthogonal linear model. 
It follows that $\sigma$ also has a 
generalized orthogonal linear model, as in Section 3. 
Indeed, by \cite[Lemma 4.2]{N10}, $\wt{\sigma} \cong \sigma^{\delta_{2m}}$, 
where 
$$\delta_{4m} = \left(\begin{array}{cccc}
       I_{2m-1} & 0 & 0 & 0\\
       0 & 0 & 1 & 0 \\
       0 & 1 & 0 & 0\\
       0 & 0 & 0 & I_{2m-1}
      \end{array}
\right),$$
and $\sigma^{\delta_{4m}}(g) = \sigma(\delta_{4m} g \delta_{4m}^{-1}) =\sigma(\delta_{4m} g \delta_{4m})$.  
So, $\sigma^{\delta_{4m}}$ has a nonzero generalized orthogonal linear model. 
Assume that $\ell$ is a nonzero generalized orthogonal linear functional of $\sigma^{\delta_{2m}}$.
Then, for any $(g_1, g_2) \in \SO_{2n} \times \SO_{4m-2n}$, 
$v \in V_{\sigma}$, 
$$\ell(\sigma(g_1, g_2) v) = \ell(\sigma^{\delta_{4m}}(g_1, g_2)^{\delta_{4m}} v) = \ell(\sigma^{\delta_{4m}}(\delta_{2n} g_1 \delta_{2n}, g_2 ) v)
= \ell(v),$$
noting that $\delta_{2n} g_1 \delta_{2n} \in SO_{2n}$. 
Hence, $\ell$ is also a nonzero generalized orthogonal linear functional of $\sigma$. Therefore, $\sigma$ has a nonzero
generalized orthogonal linear model. 

This completes the proof of Theorem \ref{main}, Part (2).

\section{About converse of Theorem \ref{main}}

In this section, we discuss some results related to the converse of Theorem \ref{main}.

First, we discuss the converse of Theorem \ref{main}, Part (1). 
Assume that $\pi$ is an irreducible admissible representation of $Sp_{4m}$ ($2m \geq n$) with a nonzero generalized symplectic linear model,
and $\sigma$ is an irreducible admissible representation of $SO_{4n}$, which is corresponding to $\pi$ under the 
local theta correspondence, i.e., 
$$\mathrm{Hom}_{SO_{4n} \times Sp_{4m}} (\omega_{\psi}, \sigma \otimes \pi) \neq 0.$$

Since $\pi$ has a nonzero generalized symplectic linear model, that is, 
$$\mathrm{Hom}_{Sp_{2n} \times Sp_{4m-2n}}(V_{\pi}, 1)  \neq 0,$$
we have that the space of $SO_{4n} \times Sp_{2n} \times Sp_{4m-2n}$-equivariant homomorphisms
\begin{equation} \label{equ11}
\mathrm{Hom}_{SO_{4n} \times Sp_{2n} \times Sp_{4m-2n}} (\omega_{\psi}, \sigma)  \neq 0.
\end{equation}

As in Section 3, let $V$ be a $4n$-dimensional vector space over $F$, with the nondegenerate symmetric from $v_{4n}$. Fix 
a basis 
$$\{e_1, \ldots, e_{2n}, e_{-2n}, \ldots, e_{-1}\}$$
 of $V$ over $F$, such that $(e_i, e_j) = (e_{-i}, e_{-j}) = 0,
(e_i, e_{-j}) = \delta_{ij}$, for $i, j = 1, \ldots, 2n$. Let $V^+ = \mathrm{Span}_F\{e_1, \ldots, e_{2n}\}$,
$V^- = \mathrm{Span}_F\{e_{-1}, \ldots, e_{-2n}\}$, then $V= V^+ + V^-$ is a complete polarization of $V$.
Fix a basis 
$$\{f_1, \ldots, f_{2m}, f_{-2m}, \ldots, f_{-1}\}$$
of $W$ over $F$, such that $(f_i, f_j) = (f_{-i}, f_{-j}) = 0,
(f_i, f_{-j}) = \delta_{ij}$, for $i, j = 1, \ldots, 2m$. 
Let 
\begin{align*}
W_1 & = \ \mathrm{Span}_F\{f_{2m-n+1}, \ldots, f_{2m}, f_{-2m}, \ldots, f_{-2m+n-1}\},\\
W_2 & = \ \mathrm{Span}_F\{f_1, \ldots, f_{2m-n}, f_{-2m+n}, \ldots, f_1\}.
\end{align*}

Then $W = W_1 + W_2$, and $\omega_{\psi} = \omega_{\psi}^1 \otimes \omega_{\psi}^2$, where
$\omega_{\psi}^1 \cong \mathcal{S}(V^- \otimes W_1)$, and $\omega_{\psi}^2 \cong \mathcal{S}(V^- \otimes W_2)$.

Then by \eqref{equ11}, 
\begin{equation} \label{equ12}
\mathrm{Hom}_{SO_{4n} \times Sp_{2n} \times Sp_{4m-2n}} 
(\omega_{\psi}^1 \otimes \omega_{\psi}^2, \sigma)  \neq 0.
\end{equation}

Let $\Theta(1_{W_1}, V) = \mathcal{S}(V^- \otimes W_1)_{Sp_{2n}}$ be the maximal quotient of $\mathcal{S}(V^- \otimes W_1)$
on which $Sp_{2n}$ acts trivially, and $\Theta(1_{W_2}, V) = \mathcal{S}(V^- \otimes W_2)_{Sp_{4m-2n}}$ be the 
maximal quotient of $\mathcal{S}(V^- \otimes W_2)$ on which $Sp_{4m-2n}$ acts trivially. Then by \eqref{equ12}, one can see that $\sigma$ must satisfy the following condition
\begin{equation*}
\mathrm{Hom}_{SO_{4n}} 
(\Theta(1_{W_1}, V) \otimes \Theta(1_{W_2}, V), \sigma)  \neq 0.
\end{equation*}

Next, we discuss the converse of Theorem \ref{main}, Part (2). 
Assume that $\sigma$ is an irreducible admissible representation of $SO_{4m}$ ($2m \geq n$) with a nonzero generalized orthogonal linear model,
and $\pi$ is an irreducible admissible representation of $Sp_{4n}$, which is corresponding to $\sigma$ under the 
local theta correspondence, that is,  
$$\mathrm{Hom}_{SO_{4m} \times Sp_{4n}} (\omega_{\psi}, \sigma \otimes \pi) \neq 0.$$

Since $\sigma$ has a generalized orthogonal linear model, that is, 
$$\mathrm{Hom}_{SO_{2n} \times SO_{4m-2n}}(V_{\sigma}, 1)  \neq 0,$$ we have that the space of $SO_{2n} \times SO_{4m-2n} \times Sp_{4n}$-equivariant homomorphisms
\begin{equation} \label{equ24}
\mathrm{Hom}_{SO_{2n} \times SO_{4m-2n} \times Sp_{4n}} (\omega_{\psi}, \pi)  \neq 0.
\end{equation}

As in Section 4, let $W$ be a $4n$-dimensional symplectic vector space over $F$, with the symplectic from $J_{4n} = \left(\begin{array}{cc}
        0 & v_{2n} \\
       -v_{2n} & 0 
      \end{array}
\right)$. Fix 
a basis 
$$\{f_1, \ldots, f_{2n}, f_{-2n}, \ldots, f_{-1}\}$$
of $W$ over $F$, such that $(f_i, f_j) = (f_{-i}, f_{-j}) = 0,
(f_i, f_{-j}) = \delta_{ij}$, for $i, j = 1, \ldots, 2n$. Let $W^+ = \mathrm{Span}_F\{f_1, \ldots, f_{2n}\}$,
$W^- = \mathrm{Span}_F\{f_{-1}, \ldots, f_{-2n}\}$, then $W= W^+ + W^-$ is a complete polarization of $W$.
Let $V$ be a $4m$-dimensional vector space over $F$, with the nondegenerate symmetric from $v_{4m}$. Fix 
a basis 
$$\{e_1, \ldots, e_{2m}, e_{-2m}, \ldots, e_{-1}\}$$
 of $V$ over $F$, such that $(e_i, e_j) = (e_{-i}, e_{-j}) = 0,
(e_i, e_{-j}) = \delta_{ij}$, for $i, j = 1, \ldots, 2m$.
Let 
\begin{align*}
V_1 & = \mathrm{Span}_F\{e_{2m-n+1}, \ldots, e_{2m}, e_{-2m}, \ldots, e_{-2m+n-1}\},\\
V_2 & = \mathrm{Span}_F\{e_1, \ldots, e_{2m-n}, e_{-2m+n}, \ldots, e_1\}.
\end{align*}
Then $V = V_1 + V_2$, and $\omega_{\psi} = \omega_{\psi}^1 \otimes \omega_{\psi}^2$, where
$\omega_{\psi}^1 \cong \mathcal{S}(V_1 \otimes W^+)$, and $\omega_{\psi}^2 \cong \mathcal{S}(V_2 \otimes W^+)$.

Then by \eqref{equ24}, 
\begin{equation} \label{equ25}
\mathrm{Hom}_{SO_{2n} \times SO_{4m-2n} \times Sp_{4n}} 
(\omega_{\psi}^1 \otimes \omega_{\psi}^2, \pi)  \neq 0.
\end{equation}

Let $\Theta(1_{V_1}, W) = \mathcal{S}(V_1 \otimes W^+)_{SO_{2n}}$ be the maximal quotient of $\mathcal{S}(V_1 \otimes W^+)$
on which $SO_{2n}$ acts trivially, and $\Theta(1_{V_2}, W) = \mathcal{S}(V_2 \otimes W^+)_{SO_{4m-2n}}$ be the 
maximal quotient of $\mathcal{S}(V_2 \otimes W^+)$ on which $SO_{4m-2n}$ acts trivially. Then by \eqref{equ25}, one can see that $\pi$ must satisfy the following condition
\begin{equation*}
\mathrm{Hom}_{Sp_{4n}} 
(\Theta(1_{V_1}, W) \otimes \Theta(1_{V_2}, W), \pi)  \neq 0.
\end{equation*}

We summarize the above discussions as follows. 

\begin{thm}\label{main2}
\begin{enumerate}
\item Assume that $\pi$ is an irreducible admissible representation of $Sp_{4m}$ $(2m \geq n)$ with a nonzero generalized symplectic linear model,
and $\sigma$ is an irreducible admissible representation of $SO_{4n}$ which is corresponding to $\pi$ under the local theta correspondence, that is,  
$$\mathrm{Hom}_{SO_{4n} \times Sp_{4m}} (\omega_{\psi}, \sigma \otimes \pi) \neq 0.$$
Then, $\sigma$ must satisfy the following condition
\begin{equation}\label{converseequ1}
\mathrm{Hom}_{SO_{4n}} 
(\Theta(1_{W_1}, V) \otimes \Theta(1_{W_2}, V), \sigma)  \neq 0.
\end{equation}
\item Assume that $\sigma$ is an irreducible admissible representation of $SO_{4m}$ $(2m \geq n)$ with a nonzero generalized orthogonal linear model,
and $\pi$ is an irreducible admissible representation of $Sp_{4n}$  which is corresponding to $\sigma$ under the 
local theta correspondence, that is,  
$$\mathrm{Hom}_{SO_{4m} \times Sp_{4n}} (\omega_{\psi}, \sigma \otimes \pi) \neq 0.$$ Then $\pi$ must satisfy the following condition
\begin{equation}\label{converseequ2}
\mathrm{Hom}_{Sp_{4n}} 
(\Theta(1_{V_1}, W) \otimes \Theta(1_{V_2}, W), \pi)  \neq 0.
\end{equation}
\end{enumerate}
\end{thm}

From Theorem \ref{main} and Theorem \ref{main2}, we have the following properties for representations of $SO_{4n}$ or $Sp_{4n}$ with nonzero generalized Shalika models. 

\begin{cor}
\begin{enumerate}
\item Assume that $\sigma$ is an irreducible admissible representation of $SO_{4n}$ with a nonzero generalized Shalika model,
and $\sigma$ occurs in the local theta correspondence of $SO_{4n}$ with $Sp_{4m}$ ($2m \geq n$). Then $\sigma$ satisfies \eqref{converseequ1}.
\item Assume that $\pi$ is an irreducible admissible representation of $Sp_{4n}$ with a nonzero generalized Shalika model,
and $\pi$ occurs in the local theta correspondence of $Sp_{4n}$ with $SO_{4m}$ ($2m \geq n$).
Then $\pi$ satisfies \eqref{converseequ2}.
\end{enumerate}
\end{cor}

\end{document}